\documentclass[11pt]{article}
\usepackage[brazil,english]{babel}
\usepackage[latin1]{inputenc}
\usepackage{fullpage}
\usepackage[dvips]{graphicx}
\usepackage{amssymb}
\usepackage{amsmath}
\usepackage{color}
\usepackage{footnote}
\usepackage{hyperref}

\newtheorem{theorem}{Theorem}[section] 
\newtheorem{corollary}[theorem]{Corollary}  
\newtheorem{proposition}[theorem]{Proposition}  
\newtheorem{lemma}[theorem]{Lemma}  
\newtheorem{definition}[theorem]{Definition} 
\newtheorem{example}[theorem]{Example}  
\newtheorem{remark}[theorem]{Remark}  

\def\beginproof{\noindent{\em Proof}.~}
\def\endproof{{\ \hfill\hbox{\fbox{}}\parfillskip 0pt}\par}

\renewcommand{\thefootnote}{\arabic{ftnote}}

\newcommand{\R}{\mathbb{R}}
\definecolor{cmg}{RGB}{250, 0, 250}

\title{On the weak stationarity conditions for Mathematical Programs with Cardinality 
Constraints: a unified approach}

\author{Evelin H. M. Krulikovski\footnotemark[3] \and Ademir A. Ribeiro\footnotemark[4]
\footnotemark[5] \and Mael Sachine\footnotemark[4]}

\begin{document}

\maketitle
\renewcommand{\thefootnote}{\fnsymbol{footnote}}
\footnotetext[3]{Graduate Program in Mathematics, Federal University of Paran\'a, 
Brazil (evelin.hmk@gmail.com).}
\footnotetext[4]{Department of Mathematics, Federal University of Paran\'a, 
Brazil (ademir.ribeiro@ufpr.br, mael@ufpr.br).}
\footnotetext[5]{Partially supported by CNPq, Brazil, Grant 309437/2016-4.}

\renewcommand{\thefootnote}{\arabic{ftnote}}

\begin{abstract}
In this paper, we study a class of optimization problems, called Mathematical Programs 
with Cardinality Constraints (MPCaC). This kind of problem is generally difficult to 
deal with, because it involves a constraint that is not continuous neither convex, 
but provides sparse solutions. Thereby we reformulate MPCaC in a suitable way, by 
modeling it as mixed-integer problem and then addressing its continuous counterpart, 
which will be referred to as relaxed problem. We investigate the relaxed problem by 
analyzing the classical constraints in two cases: linear and nonlinear. 
In the linear case, we propose a general approach and present a discussion of the 
Guignard and Abadie constraint qualifications, proving in this case that every 
minimizer of the relaxed problem satisfies the Karush-Kuhn-Tucker (KKT) conditions.
On the other hand, in the nonlinear case, we show that some standard constraint 
qualifications may be violated. Therefore, we cannot assert about KKT points. 
Motivated to find a minimizer for the MPCaC problem, we define new and weaker 
stationarity conditions, by proposing a unified approach that goes from the weakest 
to the strongest stationarity.
\end{abstract}


{\bf Subclass. }
90C30, 90C33, 90C46

\thispagestyle{plain}

\section{Introduction}
In this paper we study a class of optimization problems called {\em Mathematical Programs 
with Cardinality Constraints}, MPCaC for short, given by 
\begin{equation}
\label{prob:mpcac}
\begin{array}{cl}
\displaystyle\mathop{\rm minimize }  & f(x)  \\
{\rm subject\ to } & x \in X, \\
& \|x\|_0\leq \alpha,
\end{array}
\end{equation}
where $f:\R^n\to\R$ is a continuously differentiable function, $X\subset\R^n$ is a set 
given by equality and/or inequality constraints, $\alpha>0$ is a given natural number 
and $\|x\|_0$ denotes the cardinality of the vector $x\in\R^n$, that is, the number 
of nonzero components of $x$. Throughout this work we assume that $\alpha<n$, since 
otherwise the cardinality constraint would not have any effect. On the other side, 
if $\alpha$ is too small, the cardinality constraint may be too restrictive leading 
to an empty feasible set.

Furthermore, the main difference between problem (\ref{prob:mpcac}) and a standard 
nonlinear programming problem is that the cardinality constraint, despite of the 
notation, is not a norm, nor continuous neither convex. A classical way \cite{Bienstock} 
to deal with this difficult cardinality constraint consists of introducing binary 
variables and then rewriting the problem as a mixed-integer problem 
\begin{equation}
\label{prob:mix_int}
\begin{array}{cl}
\displaystyle\mathop{\rm minimize }_{x,y} & f(x)  \\
{\rm subject\ to } & x \in X, \\
& e^Ty\geq n-\alpha, \\
& x_iy_i=0,\; i=1,\ldots,n,\\
& y_i\in \{0,1\},\; i=1,\ldots,n,
\end{array}
\end{equation}
where $e\in\R^n$ denotes the vector of ones. 
Note that this reformulation is quite natural by noting that if a vector $x\in\R^n$ 
is such that $\|x\|_0=r\leq\alpha$, defining $y\in\R^n$ by $y_i=0$, if $x_i\neq 0$ 
and $y_i=1$, if $x_i=0$, we have $e^Ty=n-r\geq n-\alpha$ and $x_iy_i=0$ for 
all $i=1,\ldots,n$.

Alternatively to the formulation (\ref{prob:mix_int}), which is still complicated 
to deal with, in view of the binary variables, one may address its continuous 
counterpart \cite{BurdakovKanzowSchwartz16}
\begin{equation}
\label{prob:relax}
\begin{array}{cl}
\displaystyle\mathop{\rm minimize }_{x,y} & f(x)  \\
{\rm subject\ to } & x \in X, \\
& e^Ty\geq n-\alpha, \\
& x_iy_i=0,\; i=1,\ldots,n, \\
& 0\leq y_i \leq 1, \; i=1,\ldots,n,
\end{array}
\end{equation}
which will be referred to as relaxed problem and plays an important role in this 
work. As we shall see in Section \ref{sec:relations}, the 
problems (\ref{prob:mix_int}) and (\ref{prob:relax}) are closely related 
to (\ref{prob:mpcac}) in terms of feasible points and solutions. 

Despite demanding artificial variables that increase its dimension, it is precisely due to 
this augmentation that  problem (\ref{prob:relax}) has the feature of being manageable, 
in the sense of being differentiable which allows one to discuss stationarity concepts. 
This approach is common to deal with optimization 
problems \cite{Burke,GillMurraySaunders,RibeiroSachineSantos17,RibeiroSachineSantos18,Svanberg}. 

In many areas of applications of optimization we seek to find solutions with a small or a 
bounded number of nonzero components, namely sparse solutions, such as sampling signals or 
images, machine learning, subset selection in regression, portfolio problems 
\cite{BrandaBucherCervinkaSchwartz,CandesWakin,AspremontBachGhaoui,Markowitz,Miller,Tibshirani}. 
See also \cite{BeckEldar,TorrubianoMoratillaSuarez,SunZhengLi} and the references therein for 
some more ideas. 

One standard way to obtain sparse solutions consists of employing penalization 
techniques based on the $\ell_1$-norm \cite{HastieTibshiraniFriedman}. Another way is 
imposing explicitly a cardinality constraint to the problem, as the 
pioneering work \cite{Bienstock}. 
In this paper, we follow the approach that considers the cardinality constrained 
problem MPCaC, 
as \cite{BucherSchwartz,BurdakovKanzowSchwartz15,BurdakovKanzowSchwartz16}. 
Specifically, we focus on the theoretical features of the cardinality 
problem (\ref{prob:mpcac}), which may be inferred from the properties of the relaxed 
problem (\ref{prob:relax}), in view of constraint qualifications (CQ) and stationarity 
concepts. We stress that in this work we are neither concerned with applications nor with 
computational aspects or algorithmic consequences.

We consider two cases: when the set $X$ is given by linear constraints in the relaxed 
problem (\ref{prob:relax}), providing a feasible set consisting of linear (separable) 
constraints in the variables $x$ and $y$ besides the complementarity constraint; and 
when the set $X$ is given by nonlinear constraints.

In the first case, we propose a general approach that allows us to simplify the proofs 
of the results, when comparing with the ones presented in \cite{BurdakovKanzowSchwartz16}, 
as well as to prove Abadie constraint qualification (ACQ), instead of only Guignard 
constraint qualification (GCQ). We therefore conclude that 
every minimizer of the relaxed problem (\ref{prob:relax}) satisfies the KKT conditions. 
This, however, does not mean that weaker stationarity 
conditions are unnecessary or less important. Despite this is 
not the focus of this work, they are of interest from both the theoretical 
and the algorithmic viewpoint, as in the sparsity constrained 
optimization (when $X=\mathbb{R}^n$) \cite{BeckEldar,LiSong,PanXiuZhou}.

On the other hand, in the nonlinear case, we show that the most known CQs, namely LICQ 
and MFCQ, are almost never satisfied. We also prove that even a weaker condition, ACQ, 
fails to hold in a wide range of cardinality problems. Moreover, still GCQ, the weakest CQ, 
may be violated. Therefore, we cannot assert the convergence results for MPCaC in the 
same way as we usually have in the context of standard nonlinear programming, 
i.e. for KKT points.

Motivated to find a minimizer for the MPCaC problem (\ref{prob:mpcac}) for general 
constraints, we define new and weaker stationarity conditions. 
This approach is common in the literature. For example, in the 
works \cite{DussaultHaddouKadraniMigot,Outrata,ScheelScholtes} the following stationarity 
conditions are established: Weak, Clarke, Mordukhovich and Strong; while for Mathematical 
Programs with Vanishing Constraints (MPVC) the concept of $T$-stationarity was 
proposed \cite{HoheiselKanzowSchwartz}.
 
In order to find stationarity conditions for this class of problems, in the present work 
we define an auxiliary problem, namely Tightened Nonlinear Problem.  Its resulting 
formulation is similar to that made for Mathematical Programs with Complementarity 
Constraints (MPCC) \cite{Izmailov}. In this way, based on the weaker 
stationarity concepts made for the MPCC in \cite{AndreaniHaeserSecchinSilva,Ramos}, 
we propose new stationarity concepts for the class of MPCaC problems. 

Specifically, we propose a unified approach that goes from the weakest to the strongest 
stationarity for the cardinality problem with general constraints. This approach, which will 
be called $W_I$-stationarity, is based on a given set of indices $I$ such that the 
complementarity constraint is always satisfied. Moreover, different levels of stationarity 
can be obtained depending on the range for the set $I$. Besides, we prove that this condition 
is indeed weaker than the classical KKT condition, that is, every KKT point fulfils 
$W_I$-stationarity. We also point out that our definition generalizes the concepts 
of $S$- and $M$-stationarity presented in \cite{BurdakovKanzowSchwartz16} for proper 
choices of the index set $I$. In addition, specializing $W_I$-stationarity by 
considering solely the cardinality constraint, we relate our concept to the 
notions of $L$-, $N$- and $T$-stationarity discussed in \cite{BeckEldar,PanXiuZhou,LiSong} 
for sparsity constrained optimization.

We stress that although the relaxed problem (\ref{prob:relax}) resembles 
an MPCC problem, for which there is a vast literature, there are important differences 
between these two classes of optimization problems, which in turn increases the importance 
of specialized research on MPCaC problems. One of such differences is that here we only 
require positivity for one term in the complementarity constraint $x_iy_i=0$. Moreover, 
we establish results that are stronger than the corresponding ones known for MPCC's, 
as for example, the fulfillment of GCQ in the linear case (see Remark \ref{rm:MPCC} ahead).
It is also worth mentioning that, besides the usual constraint qualifications, 
also the MPCC-tailored constraint qualifications are violated for MPCaC problems in the 
general case (see \cite{BurdakovKanzowSchwartz16,CervinkaKanzowSchwartz} for a more 
detailed discussion). On the other hand, we shall conclude that $W_I$-stationarity is a 
necessary optimality condition under the MPCaC-tailored constraint qualifications proposed 
in \cite{CervinkaKanzowSchwartz}.


The paper is organized as follows: in Section \ref{sec:prelim} we establish 
some definitions, basic results and examples concerning standard nonlinear 
programming and cardinality constrained problems. In Section \ref{sec:card_linear} we 
present one of the contributions of this paper, by considering the relaxed 
problem (\ref{prob:relax}) with $X$ given by linear constraints. 
Section~\ref{sec:card_nonlinear} is devoted to our main contribution, presenting the 
analysis of the nonlinear case, including a discussion of the main constraint qualifications, 
with results, examples and counterexamples. We define weaker stationarity concepts by 
proposing a unified approach that goes from the 
weakest to the strongest stationarity. Finally, in Section \ref{sec:concl}, we discuss some 
possibilities of further research in this subject.

{\noindent\bf Notation.} Throughout this paper, we use $\|\cdot\|$ to denote the Euclidean 
norm. For vectors $x,y\in\R^n$, $x*y$ denotes the Hadamard product between $x$ and $y$, 
that is, the vector obtained by the componentwise product of $x$ and $y$. We also consider 
the following sets of indices: 
$I_{00}(x,y)=\{i\mid x_i=0,y_i=0\}$, 
$I_{\pm 0}(x,y)=\{i\mid x_i\neq 0,y_i=0\}$,
$I_{0\pm}(x,y)=\{i\mid x_i=0,y_i\neq 0\}$,
$I_{0+}(x,y)=\{i\mid x_i=0,y_i\in(0,1)\}$,
$I_{0\,>}(x,y)=\{i\mid x_i=0,y_i>0\}$,
$I_{01}(x,y)=\{i\mid x_i=0,y_i=1\}$ and 
$I_0(x)=\{i\mid x_i=0\}$. For a function $g:\R^n\to\R^m$, denote $I_g(x)=\{i\mid g_i(x)=0\}$, 
the set of active indices, and $\nabla g= (\nabla g_1,\ldots,\nabla g_m)$, the transpose of the 
Jacobian of $g$.

\section{Preliminaries}
\label{sec:prelim}
In this section we recall some basic definitions, results and examples regarding standard nonlinear 
programming (NLP) and MPCaC. Consider first the problem 
\begin{equation}
\label{prob:nlp}
\begin{array}{cl}
\displaystyle\mathop{\rm minimize } & f(x)\\
{\rm subject\ to } &  g(x)\leq 0, \\
& h(x)=0,
\end{array}
\end{equation}
where $f:\R^n\to\R$, $g:\R^n\to\R^m$ and $h:\R^n\to\R^p$ are continuously differentiable 
functions. The feasible set of the problem (\ref{prob:nlp}) is denoted by 
\begin{equation}
\label{feas_set}
\Omega=\{x\in\R^n\mid g(x)\leq 0, h(x)=0\}.
\end{equation}

%
%

\subsection{Constraint qualifications}
\label{sec:cq}
There are a lot of constraint qualifications, that is, conditions under which every minimizer 
satisfies KKT. In order to discuss some of them, let us recall the definition of cone, which 
plays an important role in this context. 

We say that a nonempty set $C\subset\R^n$ is a {\em cone} if $td\in{C}$ for all $t\geq 0$ 
and $d\in{C}$. Given a set $S\subset\R^n$, its {\em polar} is the cone 
$$
S^\circ=\{p\in\R^n\mid p^Tx\leq 0,\ \forall x\in S\}.
$$

Associated with the feasible set of the problem (\ref{prob:nlp}), we have 
the {\em tangent cone} 
$$
T_{\Omega}(\bar{x})=\left\{d\in\R^n\mid\exists(x^k)\subset\Omega\mbox{, } 
(t_k)\subset\R_+ \mbox{ : }t_k\to 0 \mbox{ and } 
\dfrac{x^k-\bar{x}}{t_k}\to d\right\}
$$
and the {\em linearized cone} 
$$
D_{\Omega}(\bar{x})=\left\{d\in\R^n\mid\nabla g_i(\bar{x})^Td \leq  0, \;
i\in I_g(\bar{x})
\mbox{ and }\nabla h(\bar{x})^Td=0 \right\}.
$$

The following basic result, whose proof is straightforward, says that we may ignore 
inactive constraints when dealing with the tangent and linearized cones. 
\begin{lemma}
\label{lm:inactive}
Consider a feasible point $\bar{x}\in\Omega$, an index set $J\supset I_g(\bar{x})$ 
and 
$$
\Omega'=\{x\in\R^n\mid g_{i}(x)\leq 0,\, i\in J, \; h(x)=0\}.
$$
Then, $T_{\Omega}(\bar{x})=T_{\Omega'}(\bar{x})$ and 
$D_{\Omega}(\bar{x})=D_{\Omega'}(\bar{x})$.
\end{lemma}

Besides the well known and strongest constraint qualifications linear 
independence constraint qualification (LICQ) and Mangasarian-Fromovitz constraint 
qualification (MFCQ), we cite here the two weakest ones.
\begin{definition}
\label{acq_gcq}
We say that Abadie constraint qualification (ACQ) holds at $\bar{x}\in\Omega$ if 
$T_{\Omega}(\bar{x})=D_{\Omega}(\bar{x})$. If 
$T_{\Omega}^\circ(\bar{x})=D_{\Omega}^\circ(\bar{x})$, we say that Guignard constraint 
qualification (GCQ) holds at $\bar{x}$.
\end{definition}

An interesting property obtained under the {\em strict complementarity} 
condition $I_{00}(\bar{x},\bar{y})=\emptyset$ says that a constraint of the form 
$x*y=0$ near 
$(\bar{x},\bar{y})$ is given only by linear constraints, as proved in the following 
result.

\begin{lemma}
\label{lm:strict_compl}
Consider the set 
$$
\Omega=\{(x,y)\in\R^n\times\R^n\mid \varphi(x,y)\leq 0,\, \rho(x,y)=0,\, x*y=0\},
$$
where $\varphi:\R^n\times\R^n\to\R^m$ and $\rho:\R^n\times\R^n\to\R^p$ are continuously 
differentiable functions. Suppose that $(\bar{x},\bar{y})\in\Omega$ satisfies 
$I_{00}(\bar{x},\bar{y})=\emptyset$ and define the set  
$$
\Omega'=\{(x,y)\in\R^n\times\R^n\mid \varphi(x,y)\leq 0,\, \rho(x,y)=0,\, 
x_{I_{0\pm}}=0,\, y_{I_{\pm 0}}=0\},
$$
where $I_{0\pm}=I_{0\pm}(\bar{x},\bar{y})$ and $I_{\pm 0}=I_{\pm 0}(\bar{x},\bar{y})$. 
Then, $T_{\Omega}(\bar{x},\bar{y})=T_{\Omega'}(\bar{x},\bar{y})$ and 
$D_{\Omega}(\bar{x},\bar{y})=D_{\Omega'}(\bar{x},\bar{y})$.
\end{lemma}
\beginproof
Note first that there exists $\delta>0$ such that 
$$
B\big((\bar{x},\bar{y}),\delta\big)\cap\Omega=B\big((\bar{x},\bar{y}),\delta\big)\cap\Omega',
$$
giving $T_{\Omega}(\bar{x},\bar{y})=T_{\Omega'}(\bar{x},\bar{y})$. 
On the other hand, denoting $\xi(x,y)=x*y$, we have  
$$
{\rm span}\left\{\nabla\xi_i(\bar{x},\bar{y}),\, i=1,\ldots,n\right\}=
{\rm span}\left\{( e_i,0),(0,e_j),\, i\in{I_{0\pm}},\,j\in{I_{\pm 0}}\right\},
$$
yielding $D_{\Omega}(\bar{x},\bar{y})=D_{\Omega'}(\bar{x},\bar{y})$.
\endproof

\subsection{Relations between the MPCaC and the reformulated problems}
\label{sec:relations}
In this section we present results that show some properties of the reformulated problems 
(\ref{prob:mix_int}) and (\ref{prob:relax}) and the equivalence between its solutions and the 
solutions of the cardinality problem (\ref{prob:mpcac}). Such results are based on 
the ones presented in \cite{BurdakovKanzowSchwartz16} and can be easily proved.

We start by noting that it is immediate that every feasible point of the mixed-integer 
problem (\ref{prob:mix_int}) is also feasible for the relaxed problem (\ref{prob:relax}), but 
the converse is clearly false. There is, however, a particular case in which the equivalence 
holds, as we can see from the next result. 

\begin{lemma}
\label{lm:card:alpha}
Let $(\bar{x},\bar{y})$ be a feasible point of the relaxed problem (\ref{prob:relax}) and 
suppose that $\|\bar{x}\|_0=\alpha$. Then, $e^T\bar{y}=n-\alpha$, $\bar{y}_i=0$ for 
$i\notin I_0(\bar{x})$ and $\bar{y}_i=1$ for $i\in I_0(\bar{x})$. 
So, $I_{00}(\bar{x},\bar{y})=\emptyset$ and, in particular, $(\bar{x},\bar{y})$ 
is feasible for the mixed-integer problem (\ref{prob:mix_int}).
\end{lemma}

Now we relate feasible points of the cardinality problem with feasible points of 
the reformulated ones.

\begin{lemma}
\label{lm:feas_problems}
Consider a point $\bar{x}\in\R^n$. 
\begin{enumerate}
\item If $\bar{x}$ is feasible for the cardinality problem (\ref{prob:mpcac}), 
then there exists $\bar{y}\in\R^n$ such that $(\bar{x},\bar{y})$ is feasible 
for the mixed-integer problem (\ref{prob:mix_int}) and, hence, feasible for 
the relaxed problem (\ref{prob:relax}). If, in addition, $\|\bar{x}\|_0=\alpha$, 
then $\bar{y}$ is unique;
\item If $(\bar{x},\bar{y})$ is feasible for (\ref{prob:relax}), then $\bar{x}$ is 
feasible for (\ref{prob:mpcac}). 
\end{enumerate}
\end{lemma}

The following theorem states that the MPCaC problem has a global minimizer if and only if 
the reformulated problems have global minimizers too. 

\begin{theorem}
\label{th:equiv1}
Consider a point $x^*\in\R^n$. 
\begin{enumerate}
\item If $x^*$ is a global solution of the problem (\ref{prob:mpcac}), there exists 
$y^*\in\R^n$ such that $(x^*,y^*)$ is a global solution of the problems 
(\ref{prob:mix_int}) and (\ref{prob:relax}). Moreover, for each reformulated 
problem, every feasible pair of the form $(x^*,\bar{y})$ is optimal;
\item If $(x^*,y^*)$ is a global solution of (\ref{prob:mix_int}) or (\ref{prob:relax}), 
$x^*$ is a global solution of~(\ref{prob:mpcac}). 
\end{enumerate}
\end{theorem}

As a consequence of Theorem \ref{th:equiv1}, we have that every global solution 
of (\ref{prob:mix_int}) is also a global solution of  (\ref{prob:relax}). 
However, the converse is not necessarily true, as we can see in the example below. 

\begin{example}
\label{ex:ynotbinary}
Consider the relaxed problem 
$$
\begin{array}{cl}
\displaystyle\mathop{\rm minimize }_{x,y\in\R^3} & (x_1-1)^2+(x_2-1)^2+x_3^2 \\
{\rm subject\ to } & x_1\leq 0, \\
& y_1+y_2+y_3\geq 1, \\
& x_iy_i=0,\; i=1,2,3, \\
& 0\leq y_i \leq 1, \; i=1,2,3.
\end{array}
$$
Given any $t\in[0,1]$, the pair $(x^*,y^*)$, with $x^*=(0,1,0)$ and $y^*=(1-t,0,t)$, 
is a global solution of the relaxed problem, but for $t\in(0,1)$ this point is not 
even feasible for the mixed-integer problem (\ref{prob:mix_int}). 
\end{example}

Now, let us discuss an existence result for the MPCaC and the reformulated problems. 
For this purpose, note first the closedness of the set defined by the cardinality 
constraint. Indeed, despite the fact that the function $x\mapsto\|x\|_0$
is not continuous, it is lower semicontinuous and hence, given $\alpha\geq 0$, the 
level set $\mathcal{C}=\{x\in\R^n\mid\|x\|_0\leq \alpha\}$ is 
closed (see \cite[Example 2.11]{Beck}). Therefore, another consequence of 
Theorem \ref{th:equiv1} is the following result. 

\begin{theorem}\cite{BurdakovKanzowSchwartz16}
\label{th:existence} 
Suppose that the feasible set $\Omega_0=\{x\in X\mid\|x\|_0\leq\alpha\}$ of the 
cardinality-constrained problem (\ref{prob:mpcac}) is nonempty and that $X$ is compact. Then 
the problems (\ref{prob:mpcac}), (\ref{prob:mix_int}) and (\ref{prob:relax}) have a nonempty 
solution set.  
\end{theorem}

Now, let us analyze the relations among the problems when we consider local solutions. We 
shall see that, differently from the global case, part of the equivalence is lost, but the 
relation from MPCaC to the relaxed problem remains valid. 

\begin{theorem}
\cite{BurdakovKanzowSchwartz16}
\label{th:local}
Let $x^*\in\R^n$ be a local minimizer of (\ref{prob:mpcac}). Then there exists a vector 
$y^*\in\R^n$ such that the pair $(x^*,y^*)$ is a local minimizer of (\ref{prob:relax}).
\end{theorem}

The next example shows that the converse of the above result is not valid. 

\begin{example}
\label{ex:local:noteq}
Consider the MPCaC and the corresponding relaxed problem 
$$
\begin{array}{lr}
\begin{array}{cl}
\displaystyle\mathop{\rm minimize }_{x\in\R^3} & x_1^2+(x_2-1)^2+(x_3-1)^2  \\
{\rm subject\ to } & x_1\leq 0, \\
& \|x\|_0\leq 2, \\ & \\ &
\end{array}
&
\hspace{.2cm}
\begin{array}{cl}
\displaystyle\mathop{\rm minimize }_{x,y\in\R^3} & x_1^2+(x_2-1)^2+(x_3-1)^2 \\
{\rm subject\ to } & x_1\leq 0, \\
& y_1+y_2+y_3\geq 1, \\
& x_iy_i=0,\; i=1,2,3, \\
& 0\leq y_i \leq 1, \; i=1,2,3.
\end{array}
\end{array}
$$
Fix $t\in(0,1)$ and define $x^*=(0,1,0)$ and $y^*=(1-t,0,t)$. We claim that 
the pair $(x^*,y^*)$ is a local solution of the relaxed problem. Indeed, if $(x,y)$ is 
sufficiently close to $(x^*,y^*)$, then $y_1\neq 0$ and $y_3\neq 0$, which implies that 
$x_1=0$ and $x_3=0$. So, the objective value at $(x,y)$ is $(x_2-1)^2+1\geq 1$, proving the 
claim. Nevertheless, $x^*$ is not a local minimizer of the MPCaC, because we can consider a 
point $x_{\delta}=(0,1,\delta)$ as close to $x^*=(0,1,0)$ as we want whose objective value 
at this point is $(\delta-1)^2<1$. 
\end{example}

Note that, in the above example, there are infinitely many vectors $y^*$ such that 
$(x^*,y^*)$ is a local solution of the relaxed problem. This is the reason for $x^*$ not 
to be a local minimizer of the MPCaC, as we can see from the next result. 
\begin{theorem}
\label{th:local:converse}
Let $(x^*,y^*)$ be a local minimizer of problem (\ref{prob:relax}). Then $\|x^*\|_0=\alpha$ 
if and only if $y^*$ is unique, that is, if there is exactly one $y^*$ such that 
$(x^*,y^*)$ is a local minimizer of (\ref{prob:relax}). In this case, the components of 
$y^*$ are binary and  $x^*$ is a local minimizer of (\ref{prob:mpcac}).
\end{theorem}
\beginproof
The ``only if'' part and the claim that $y^*$ is a binary vector follow directly from 
Lemma \ref{lm:card:alpha}. The ``if'' part and the proof that $x^*$ is a local minimizer of 
(\ref{prob:mpcac}) is given in \cite{BurdakovKanzowSchwartz16}.
\endproof 

\section{MPCaC: the Linear Case}
\label{sec:card_linear}
It is well known that a set defined by linear constraints naturally satisfies a constraint 
qualification. In particular, there holds ACQ for this kind of constraints. Now we discuss 
what happens if we consider the relaxed problem (\ref{prob:relax}) with $X$ given by linear 
constraints. Note that in this case, we have linear (separable) constraints together with 
a coupling complementarity constraint. 

In fact, we present in this section one of the contributions of this paper, establishing more 
general results from which we can derive properties for the specific problem (\ref{prob:relax}). 
Moreover, this general approach enables us to simplify the proofs, when comparing with 
the ones presented in \cite{BurdakovKanzowSchwartz16}, as well as to discuss also ACQ, 
instead of only GCQ. The major difference between our approach and the strategy used 
in \cite{BurdakovKanzowSchwartz16} for proving Guignard is that they use partitions of 
the index set $I_{00}$ to construct decompositions of the feasible set and the corresponding 
cones in terms of simpler sets, whereas we decompose an arbitrary vector of the linearized 
cone as a sum of two vectors belonging to the tangent cone, making the proof very simple. 
Furthermore, we provide here also the analysis of Abadie condition under the strict 
complementarity condition. 

For this purpose, consider the feasible set 
\begin{equation}
\label{feas_set_linear}
\Omega=\{(x,y)\in\R^n\times\R^n\mid Ax=b,\; \tilde{A}x\leq\tilde{b},\; 
My=r,\; \tilde{M}y\leq\tilde{r},\; x*y=0\},
\end{equation}
where the matrices $A$, $\tilde{A}$, $M$ and $\tilde{M}$ and the vectors $b$, $\tilde{b}$, 
$r$ and $\tilde{r}$ have appropriate dimensions.

We start by giving sufficient conditions for ACQ to be satisfied.
\begin{theorem}
\label{th:lACQ}
Consider the set $\Omega$, defined in (\ref{feas_set_linear}), and a point 
$(\bar{x},\bar{y})\in\Omega$. If $I_{00}(\bar{x},\bar{y})=\emptyset$, 
then ACQ holds at $(\bar{x},\bar{y})$.
\end{theorem}
\beginproof
Follows directly from Lemma \ref{lm:strict_compl}.
\endproof

\begin{remark}
\label{rm:ACQ_I00}
Note that the validity of Theorem \ref{th:lACQ} does not depend on the 
separability of the linear constraints, that is, it is valid for the more general set 
\begin{equation}
\label{feas_set_linear_general}
\{(x,y)\in\R^n\times\R^n\mid Bx+Cy=c,\; \bar{B}x+\bar{C}y\leq\bar{c},\; x*y=0\},
\end{equation}
where the matrices $B$, $C$, $\bar{B}$ and $\bar{C}$ and the vectors $c$ and $\bar{c}$ 
have appropriate dimensions. On the other hand, we 
have nothing to say if $I_{00}(\bar{x},\bar{y})\neq\emptyset$. For example, 
letting $(\bar{x},\bar{y})=(0,0)$, ACQ is satisfied for 
$$
\Omega=\{(x,y)\in\R^2\mid x=0,\; 0\leq y\leq 1,\; xy=0\}.
$$ 
On the other hand, if 
$$
\Omega=\{(x,y)\in\R^2\mid 0\leq x\leq 1,\; 0\leq y\leq 1,\; xy=0\},
$$ 
ACQ does not hold. 
\end{remark}

Next result is more precise and tell us everything if we replace ACQ by GCQ. 

\begin{theorem}
\label{th:lGCQ}
Consider the set $\Omega$, defined in (\ref{feas_set_linear}).
Then, every feasible point $(\bar{x},\bar{y})\in\Omega$ satisfies GCQ. 
\end{theorem}
\beginproof
We invoke Lemma \ref{lm:inactive} to assume without loss of generality that there is 
no inactive constraint at $(\bar{x},\bar{y})$. Denote $\rho(x,y)=Ax-b$, 
$\tilde\rho(x,y)=\tilde{A}x-\tilde{b}$, $\zeta(x,y)=My-r$, 
$\tilde\zeta(x,y)=\tilde{M}x-\tilde{r}$ and $\xi(x,y)=x*y$. When there is no chance for 
ambiguity, we sometimes suppress the argument and, for example, write $I_{00}$ for 
$I_{00}(\bar{x},\bar{y})$, $I_{0\pm}$ for $I_{0\pm}(\bar{x},\bar{y})$ and so on.
Consider then $d=(u,v)\in D_{\Omega}(\bar{x},\bar{y})$ arbitrary. We claim that the 
vectors $(u,0)$ and $(0,v)$ belong to $T_{\Omega}(\bar{x},\bar{y})$. 
Indeed, since 
$
\bar{y}_iu_i+\bar{x}_iv_i=\nabla\xi_i(\bar{x},\bar{y})^Td=0
$
for all $i=1,\ldots,n$, we have $u_{I_{0\pm}}=0$ and $v_{I_{\pm0}}=0$. Thus, 
the sequences $t_k={1}/{k}$ and $(x^k,y^k)=(\bar{x}+t_ku,\bar{y})$ satisfy 
$$
x^k_{I_{0\pm}}=0\,, \;\; y^k_{I_{\pm0}\cup I_{00}}=0\quad\mbox{and}\quad
\dfrac{(x^k,y^k)-(\bar{x},\bar{y})}{t_k}\to(u,0).
$$
Moreover, since 
$(u,v)\in D_{\Omega}(\bar{x},\bar{y})$, we obtain 
$$
Au=0,\;\;\tilde{A}u\leq 0,\;\;Mv=0\quad\mbox{and}\quad \tilde{M}v\leq 0
$$
implying that $Ax^k=b$\,, $\tilde{A}x^k\leq\tilde{b}$\,, $My^k=r$ and 
$\tilde{M}y^k\leq\tilde{r}$. 
So, $(x^k,y^k)\subset\Omega$ and then $(u,0)\in T_{\Omega}(\bar{x},\bar{y})$. 
The fact that $(0,v)\in T_{\Omega}(\bar{x},\bar{y})$ can be proved analogously.

Finally, to establish the relation 
$T_{\Omega}^\circ(\bar{x},\bar{y})=D_{\Omega}^\circ(\bar{x},\bar{y})$ 
consider $p\in T_{\Omega}^\circ(\bar{x},\bar{y})$ and $d\in D_{\Omega}(\bar{x},\bar{y})$ 
arbitrary. As seen above, we can write $d=d^1+d^2$, with 
$d^1,d^2\in T_{\Omega}(\bar{x},\bar{y})$. Thus, $p^Td=p^Td^1+p^Td^2\leq 0$.
\endproof

It should be noted that, contrary to what occurs in Theorem \ref{th:lACQ}, the above 
result cannot be generalized for the set defined in (\ref{feas_set_linear_general}), 
as can be seen in the following example. 
\begin{example}
\label{ex_linearnotGuignard}
Consider the set 
$$
\Omega=\{(x,y)\in\R^2\mid x\geq 0,\; y\geq 0,\; -x+y\leq 0,\; xy=0\}.
$$ 
It can be seen that 
$$
T_{\Omega}(0)=\left\{d\in\R^2\mid d_1\geq 0,\, d_2=0\right\}\quad\mbox{and}\quad 
D_{\Omega}(0)=\left\{d\in\R^2\mid 0\leq d_2\leq d_1 \right\}.
$$
Hence, $T_{\Omega}^\circ(0)\neq D_{\Omega}^\circ(0)$.
\end{example}

\begin{remark}
\label{rm:MPCC}
The above example, incidentally, points out a significant difference between the 
class of problems we are considering in this paper, MPCaC, and the closely related 
problems MPCC's. Example \ref{ex_linearnotGuignard} is an MPCC problem and, although 
all but the last constraint are linear, GCQ does not hold. On the other hand, as 
established in Theorem \ref{th:GCQ} below, every feasible point of an MPCaC, with 
$X$ defined by linear constraints, satisfies Guignard.
\end{remark}

Now we obtain, as a direct consequence of the above results, the constraint qualification 
analysis for the relaxed problem (\ref{prob:relax}) in the linear case.
\begin{theorem}
\label{th:GCQ}
Consider the problem (\ref{prob:relax}), with $X$ defined by linear 
(equality and/or inequality) constraints, and its feasible set 
$$
\Omega=\{(x,y)\in X\times\R^n\mid e^Ty\geq n-\alpha,\; x*y=0,\; 0\leq y\leq e\}.
$$
Then, every feasible point $(\bar{x},\bar{y})\in\Omega$ satisfies GCQ. Moreover, If 
$I_{00}(\bar{x},\bar{y})=\emptyset$, then ACQ holds at $(\bar{x},\bar{y})$.
\end{theorem}
\beginproof
Follows directly from Theorems \ref{th:lACQ} and \ref{th:lGCQ}.
\endproof
\begin{corollary}
\label{c:ACQ}
In the context of Theorem \ref{th:GCQ}, if $\|\bar x\|_0=\alpha$, then 
ACQ holds at $(\bar{x},\bar{y})$.
\end{corollary}
\beginproof
Follows directly from Lemma \ref{lm:card:alpha}.
\endproof

As we pointed out in Remark \ref{rm:ACQ_I00}, the condition $I_{00}(\bar{x},\bar{y})=\emptyset$ 
is sufficient but not necessary for ACQ. There is, however, a situation in which the equivalence 
holds.
\begin{proposition}
\label{prop:ACQ}
Consider the problem (\ref{prob:relax}) with $X=\R^n$ and its feasible set 
$$
\Omega=\{(x,y)\in\R^n\times\R^n\mid e^Ty\geq n-\alpha,\; x*y=0,\; 0\leq y\leq e\}.
$$
Given an arbitrary feasible point $(\bar{x},\bar{y})\in\Omega$, ACQ holds at 
$(\bar{x},\bar{y})$ if and only if $I_{00}(\bar{x},\bar{y})=\emptyset$.
\end{proposition}
\beginproof
Denote $\theta(x,y)=n-\alpha-e^Ty$, $H(x,y)=-y$, $\tilde H(x,y)=y-e$ and $\xi(x,y)=x*y$. 
Then, given $d=(u,v)\in\R^n\times\R^n$, we have 
\begin{subequations}
\begin{align}
\nabla\theta(\bar{x},\bar{y})^Td=-e^Tv, \label{gradtheta} \\ 
\nabla H_i(\bar{x},\bar{y})^Td=-v_i, \label{gradH} \\ 
\nabla\tilde H_i(\bar{x},\bar{y})^Td=v_i, \label{gradHtil} \\ 
\nabla\xi_i(\bar{x},\bar{y})^Td=\bar{y}_iu_i+\bar{x}_iv_i \label{gradxii} .
\end{align}
\end{subequations}
Assume first that $I_{00}(\bar{x},\bar{y})\neq\emptyset$, take an index 
$\ell\in I_{00}(\bar{x},\bar{y})$ and consider the vector 
$\bar{d}=(e_\ell,e_\ell)\in\R^n\times\R^n$. Let us prove that 
$\bar{d}\in D_{\Omega}(\bar{x},\bar{y})$. We have 
$$
\nabla\theta(\bar{x},\bar{y})^T\bar{d}=-1\quad\mbox{and}\quad 
\nabla H_i(\bar{x},\bar{y})^T\bar{d}\leq 0
$$
for all $i$. Moreover, the constraint $\tilde H_\ell$ is inactive at $(\bar{x},\bar{y})$ and 
$\nabla\tilde H_i(\bar{x},\bar{y})^T\bar{d}=0$ for all $i\neq\ell$. Note also that 
$\nabla\xi_i(\bar{x},\bar{y})^T\bar{d}=0$ for all $i$. Therefore, 
$\bar{d}\in D_{\Omega}(\bar{x},\bar{y})$. We claim that 
$\bar{d}\notin T_{\Omega}(\bar{x},\bar{y})$. Indeed, given any 
$d=(u,v)\in T_{\Omega}(\bar{x},\bar{y})$, there exist sequences 
$(x^k,y^k)\subset\Omega$ and $t_k\to 0$ such that 
$\dfrac{(x^k,y^k)-(\bar{x},\bar{y})}{t_k}\to(u,v)$. This implies that 
$$
\dfrac{x_\ell^k}{t_k}=\dfrac{x_\ell^k-\bar{x}_\ell}{t_k}\to u_\ell\quad\mbox{and}\quad
\dfrac{y_\ell^k}{t_k}=\dfrac{y_\ell^k-\bar{y}_\ell}{t_k}\to v_\ell.
$$
So, $0=\dfrac{x_\ell^ky_\ell^k}{t_k^2}\to u_\ell v_\ell$, yielding $u_\ell v_\ell=0$. 
Therefore, $\bar{d}=(e_\ell,e_\ell)\notin T_{\Omega}(\bar{x},\bar{y})$ and, hence, ACQ does 
not hold at $(\bar{x},\bar{y})$. 
The converse follows directly from Theorem \ref{th:lACQ}.
\endproof

A concluding remark of this section is that, in view of Theorem \ref{th:GCQ}, with $X$ 
defined by linear constraints, every minimizer of the relaxed problem (\ref{prob:relax}) 
satisfies the KKT conditions. This fact, however, does not mean that weaker stationarity 
conditions are unnecessary or less important. They are of interest from both the theoretical 
and the practical viewpoint, as in the sparsity constrained optimization (when there is only 
the cardinality constraint). See \cite{BeckEldar,PanXiuZhou,LiSong} and references therein 
for a more detailed discussion. 

We now turn our attention to the general nonlinear case, to be discussed in the next section. 

\section{MPCaC: the Nonlinear Case}
\label{sec:card_nonlinear}
In this section we present the main contribution of this paper. We propose a unified 
approach that goes from the weakest to the strongest stationarity for the cardinality problem 
with general constraints. This approach, which will be called $W_I$-stationarity, is based 
on a given set of indices $I$ such that the 
complementarity constraint is always satisfied. Moreover, different levels of stationarity 
can be obtained depending on the range for the set $I$. Besides, we prove that this condition 
is indeed weaker than the classical KKT condition, that is, every KKT point fulfils 
$W_I$-stationarity. We also point out that our definition generalizes the concepts 
of $S$- and $M$-stationarity presented in \cite{BurdakovKanzowSchwartz16} for a 
proper choice of the index set $I$.

For this purpose, consider the MPCaC problem (\ref{prob:mpcac}) with 
\begin{equation}
\label{setX}
X=\{x\in\R^n\mid g(x)\leq 0, h(x)=0\},
\end{equation}
where $g:\R^n\to\R^m$ and $h:\R^n\to\R^p$ are continuously differentiable functions. 

Differently from the linear case, when the set $X$ is given by nonlinear constraints 
we cannot guarantee that the standard constraints qualifications are satisfied for the 
relaxed problem (\ref{prob:relax}), as discussed in Section \ref{sec:CQ} below, 
making convenient the study of weaker stationarity conditions, to be addressed in 
Section \ref{sec:stationarity}.

\subsection{Constraint qualifications for MPCaC}
\label{sec:CQ}
Here we show that the most known constraint qualifications, LICQ 
and MFCQ, are not satisfied almost everywhere. We also prove that even a weaker 
condition, ACQ, fails to hold in a wide range of cardinality problems. 
Moreover, still GCQ, the weakest constraint qualification, may be violated. 
This issues are related to the problematic constraints $x_iy_i=0$, $i=1,\ldots,n$. 

\begin{proposition}
\label{prop:notMFCQ}
Let $(\bar{x},\bar{y})$ be a feasible point of the problem (\ref{prob:relax}) and suppose that 
$\bar{x}_\ell\neq 0$ for some index $\ell\in\{1,\ldots,n\}$. Then $(\bar{x},\bar{y})$ does not 
satisfy MFCQ and, therefore, it does not satisfy LICQ.	
\end{proposition}
\beginproof
Denote $\xi_i(x,y)=x_iy_i$ and $H_i(x,y)=-y_i$. Given 
$d=(u,v)\in\R^n\times\R^n$, we have 
$\nabla\xi_i(\bar{x},\bar{y})^Td=\bar{y}_iu_i+\bar{x}_iv_i$ and 
$\nabla H_i(\bar{x},\bar{y})^Td=-v_i$. Since $\bar{x}_\ell\neq 0$, there holds 
$\bar{y}_\ell=0$, which implies that the constraint $H_\ell$ is active at 
$(\bar{x},\bar{y})$ and   
$$
\nabla\xi_\ell(\bar{x},\bar{y})^Td=\bar{x}_\ell v_\ell\quad\mbox{and}\quad 
\nabla H_\ell(\bar{x},\bar{y})^Td=-v_\ell.
$$
So, there is no $d\in\R^n\times\R^n$ satisfying  
$\nabla\xi_\ell(\bar{x},\bar{y})^Td=0$ and $\nabla H_\ell(\bar{x},\bar{y})^Td<0$ at the 
same time. This means that $(\bar{x},\bar{y})$ cannot satisfy MFCQ and, hence, it does not 
satisfy LICQ as well.
\endproof

Even if we consider weaker conditions than LICQ and MFCQ, they may fail to hold, as we saw 
in Proposition \ref{prop:ACQ}. In that situation, despite considering the simplest MPCaC problem, 
without constraints other than the cardinality constraint itself, ACQ does not hold if there 
is an index $i$ for which $\bar{x}_i=\bar{y}_i=0$. 

We can go further and see that even GCQ, the weakest constraint qualification, may not be 
fulfilled.

\begin{example}
\label{ex:notGCQ}
Consider the MPCaC and the corresponding relaxed problem 
$$
\begin{array}{lr}
\begin{array}{cl}
\displaystyle\mathop{\rm minimize }_{x\in\R^2} & x_1+x_2  \\
{\rm subject\ to } & -x_1+x_2^2\leq 0, \\
& \|x\|_0\leq 1, \\ & \\ &
\end{array}
&
\hspace{.5cm}
\begin{array}{cl}
\displaystyle\mathop{\rm minimize }_{x,y\in\R^2} & x_1+x_2 \\
{\rm subject\ to } & -x_1+x_2^2\leq 0, \\
& y_1+y_2\geq 1, \\
& x_iy_i=0,\; i=1,2, \\
& 0\leq y_i \leq 1, \; i=1,2.
\end{array}
\end{array}
$$
Note that $x^*=(0,0)$ is the unique global solution of the cardinality problem and, defining 
$y^*=(1,0)$, the pair $(x^*,y^*)$ is a global solution of the relaxed problem. However, this 
pair does not satisfy GCQ, because it is not KKT. 
\end{example}

The violation of GCQ in the previous example is due to the existence of a nonlinear constraint, 
since in the linear case GCQ is always satisfied. Moreover, the example shows that the classical 
stationarity conditions may not be able to detect the solution. 

\subsection{Weak stationarity conditions for MPCaC}
\label{sec:stationarity}
As we have seen before, except in special cases, e.g., where $X$ is polyhedral convex, we do not 
have a constraint qualification for the relaxed problem (\ref{prob:relax}). So, the standard KKT 
conditions are not necessary optimality conditions. 

Thus, in this section we define weaker stationarity concepts to deal with this class of problems. 
In fact, we propose a unified approach that goes from the weakest to the strongest stationarity.

For ease of presentation we consider the functions (some of which already seen in the 
proof of Proposition \ref{prop:ACQ}) $\theta:\R^n\times\R^n\to\R$, 
$H,\tilde{H},G:\R^n\times\R^n\to\R^n$ given by $\theta(x,y)=n-\alpha-e^Ty$, $H(x,y)=-y$, 
$\tilde{H}(x,y)=y-e$ and $G(x,y)=x$. Then we can rewrite the relaxed problem (\ref{prob:relax}) 
as 
\begin{equation}
\label{prob:relax1}
\begin{array}{cl}
\displaystyle\mathop{\rm minimize }_{x,y} & f(x)  \\
{\rm subject\ to } & g(x)\leq 0, h(x)=0, \\
& \theta(x,y)\leq 0, \\
& H(x,y)\leq 0, \tilde{H}(x,y)\leq 0, \\
& G(x,y)*H(x,y)=0.
\end{array}
\end{equation}

Given a feasible point $(\bar{x},\bar{y})$ for the problem (\ref{prob:relax1}) and a set of 
indices $I$ such that 
\begin{equation}
\label{indexI}
I_{0+}(\bar{x},\bar{y})\cup I_{01}(\bar{x},\bar{y})\subset I\subset I_0(\bar{x}),
\end{equation}
we have that $i\in I$ or $i\in I_{00}(\bar{x},\bar{y})\cup I_{\pm 0}(\bar{x},\bar{y})$
for all $i\in\{1,\ldots,n\}$. Thus, $G_i(\bar{x},\bar{y})=0$ or $H_i(\bar{x},\bar{y})=0$. 
This suggests to consider an auxiliary problem by removing the problematic constraint 
$G(x,y)*H(x,y)=0$ and including other ones that ensure the null product. We then define 
the $I$-{\em Tightened} Nonlinear Problem at $(\bar{x},\bar{y})$ by
\begin{equation}
\label{prob:tight}
\begin{array}{cl}
\displaystyle\mathop{\rm minimize }_{x,y} & f(x)  \\
{\rm subject\ to } & g(x)\leq 0, h(x)=0, \\
& \theta(x,y)\leq 0, \\
& \tilde{H}(x,y)\leq 0, \\
& H_i(x,y)\leq 0,\; i\in I_{0+}(\bar{x},\bar{y})\cup I_{01}(\bar{x},\bar{y}), \\
& H_i(x,y)=0,\; i\in I_{00}(\bar{x},\bar{y})\cup I_{\pm0}(\bar{x},\bar{y}), \\
& G_i(x,y)=0,\; i \in I.
\end{array}
\end{equation}
This problem will be also indicated by TNLP$_{I}(\bar{x},\bar{y})$ 
and, when there is no chance for ambiguity, it will be referred simply to 
as {\em tightened problem}. Note that we tighten only those constraints that are 
involved in the complementarity constraint $G(x,y)*H(x,y)=0$, by 
converting the active inequalities $H_i$'s into equalities and incorporating the equality 
constraints $G_i$'s. The upper bound $I_0(\bar{x})$ for the range of $I$ guarantees 
that we do not incorporate a constraint $G_i(x,y)=0$ for some $i$ such that 
$G_i(\bar{x},\bar{y})=\bar{x}_i\neq 0$.

The following lemma is a straightforward consequence of the definition of 
TNLP$_{I}(\bar{x},\bar{y})$. 

\begin{lemma}
\label{lm:tnlp1}
Consider the tightened problem (\ref{prob:tight}). Then, 
\begin{enumerate}
\item\label{inactiveH} the inequalities defined by $H_i$, 
$i\in I_{0+}(\bar{x},\bar{y})\cup I_{01}(\bar{x},\bar{y})$, are inactive at $(\bar{x},\bar{y})$;
\item\label{xybar_tnlp} $(\bar{x},\bar{y})$ is feasible for TLNP$_{I}(\bar{x},\bar{y})$; 
\item every feasible point of (\ref{prob:tight}) is feasible for (\ref{prob:relax1});
\item if $(\bar{x},\bar{y})$ is a global (local) minimizer of (\ref{prob:relax1}), then 
it is also a global (local) minimizer of TNLP$_{I}(\bar{x},\bar{y})$.
\end{enumerate}
\end{lemma}

The Lagrangian function associated with TNLP$_{I}(\bar{x},\bar{y})$ is the function 
$$
\mathcal{L}_{I}:\R^n\times\R^n\times\R^m\times\R^p\times\R\times\R^n
\times\R^n\times\R^{|I|}\to\R
$$ 
given by 
\begin{align*}
\mathcal{L}_{I}(x,y,\lambda^g,\lambda^h,\lambda^{\theta},\lambda^H,
\lambda^{\tilde{H}},\lambda_I^G)=f(x)+(\lambda^g)^Tg(x)+(\lambda^h)^Th(x)+
\lambda^{\theta}\theta(x,y) \\ 
+(\lambda^H)^TH(x,y)+(\lambda^{\tilde{H}})^T\tilde{H}(x,y)+(\lambda_I^G)^TG_{I}(x,y).
\end{align*}

Note that the tightened problem, and hence its Lagrangian, depends on the index set $I$, which 
in turn depends on the point $(\bar{x},\bar{y})$. It should be also noted that 
\begin{align*}
\nabla_{x,y}\mathcal{L}_{I}({x},{y},\lambda)=
\left(\begin{array}{c}\nabla f({x}) \\ 0 \end{array}\right)+
\sum_{i=1}^{m}\lambda_i^g\left(\begin{array}{c}\nabla g_i({x}) \\ 0
\end{array}\right)+\sum_{i=1}^{p}\lambda_i^h
\left(\begin{array}{c}\nabla h_i({x}) \\ 0 \end{array}\right)  \\
+\lambda^{\theta}\nabla\theta({x},{y})
+\sum_{i=1}^{n}\lambda_i^H \nabla H_i({x},{y})
+\sum_{i=1}^{n}\lambda_i^{\tilde{H}}\nabla \tilde{H}_i({x},{y}) \\
+\sum_{i\in I}\lambda_i^G\nabla G_i({x},{y})
=\left(\begin{array}{c} u \\ v \end{array}\right),
\end{align*}
with 
\begin{align}
u=\nabla f({x})+\sum_{i=1}^{m}\lambda_i^g\nabla g_i({x})+
\sum_{i=1}^{p}\lambda_i^h\nabla h_i({x})+\sum_{i\in I}\lambda_i^Ge_i \label{nablaLx} \\ 
v=-\lambda^{\theta}e-\lambda^H+\lambda^{\tilde{H}}. \label{nablaLy} 
\end{align}

\subsubsection{Weak stationarity}
\label{sec:weak}
Our weaker stationarity concept for the relaxed problem (\ref{prob:relax1}) is then defined in 
terms of the tightened problem as follows. 
\begin{definition}
\label{def:wstat_xy}
Consider a feasible point $(\bar{x},\bar{y})$ of the relaxed problem (\ref{prob:relax1}) 
and a set of indices $I$ satisfying (\ref{indexI}). We say that $(\bar{x},\bar{y})$ 
is $I$-weakly stationary ($W_{I}$-stationary) for this problem if there exists a vector 
$$
\lambda=(\lambda^g,\lambda^h,\lambda^{\theta},\lambda^H,\lambda^{\tilde{H}},\lambda^G)\in
\R_+^m\times\R^p\times\R_+\times\R^n\times\R_+^n\times\R^{|I|}
$$ 
such that
\begin{enumerate}
\item\label{lagrangian0} $\nabla_{x,y}{\cal{L}}_{I}(\bar{x},\bar{y},\lambda)=0$;
\item\label{g0} $(\lambda^g)^Tg(\bar{x})=0$;
\item\label{theta0} $\lambda^{\theta}\theta(\bar{x},\bar{y})=0$;
\item\label{htilde0} $(\lambda^{\tilde{H}})^T\tilde{H}(\bar{x},\bar{y})=0$;
\item\label{h0} $\lambda^H_i=0$ for all $i\in I_{0+}(\bar{x},\bar{y})\cup I_{01}(\bar{x},\bar{y})$.
\end{enumerate}    
\end{definition}
\begin{remark}
\label{rm:wstat}
In view of (\ref{nablaLx}) and (\ref{nablaLy}), the first item of Definition \ref{def:wstat_xy} 
means that 
\begin{align}
\nabla f(\bar{x})+\sum_{i=1}^{m}\lambda_i^g\nabla g_i(\bar{x})+
\sum_{i=1}^{p}\lambda_i^h\nabla h_i(\bar{x})+\sum_{i\in I}\lambda_i^Ge_i=0 
\label{nablaLx0} \\ -\lambda^{\theta}e-\lambda^H+\lambda^{\tilde{H}}=0. \label{nablaLy0} 
\end{align}
Items (\ref{g0}), (\ref{theta0}) and (\ref{htilde0}) represent the standard KKT complementarity 
conditions for the inequality constraints $g(x)\leq 0$, $\theta(x,y)\leq 0$ and 
$\tilde{H}(x,y)\leq 0$, respectively, of the tightened problem (\ref{prob:tight}). 
In view of Lemma \ref{lm:tnlp1}(\ref{inactiveH}), the last item also represents KKT 
complementarity conditions for the constraints $H_i(x,y)\leq 0$, 
$i\in I_{0+}(\bar{x},\bar{y})\cup I_{01}(\bar{x},\bar{y})$, of the tightened problem.
\end{remark}

As an immediate consequence of Remark \ref{rm:wstat} we have the following characterization of 
$W_{I}$-stationarity for the relaxed problem in terms of stationarity for the tightened problem.
\begin{proposition}
\label{prop:wstat_kkttnlp}
Let $(\bar{x},\bar{y})$ be feasible for the relaxed problem (\ref{prob:relax1}). Then, 
$(\bar{x},\bar{y})$ is $W_{I}$-stationary if and only if it is a KKT point for the tightened 
problem (\ref{prob:tight}).
\end{proposition}
\beginproof
Follows from the feasibility of $(\bar{x},\bar{y})$, stated in 
Lemma \ref{lm:tnlp1}(\ref{xybar_tnlp}), the comments in Remark \ref{rm:wstat}, the fact that 
$\mathcal{L}_{I}$ is the Lagrangian of TNLP$_{I}(\bar{x},\bar{y})$ and that the 
nonnegativeness of the multipliers corresponding to the inequality constraints $H_i(x,y)\leq 0$, 
$i\in I_{0+}(\bar{x},\bar{y})\cup I_{01}(\bar{x},\bar{y})$, is equivalent to the last item of 
Definition~\ref{def:wstat_xy}, because of Lemma \ref{lm:tnlp1}(\ref{inactiveH}).
\endproof

Note that in view of Proposition \ref{prop:wstat_kkttnlp} we could have defined 
$W_{I}$-stationarity simply as KKT for the tightened problem (\ref{prob:tight}). Nevertheless, 
we prefer as in Definition \ref{def:wstat_xy} in order to have its last condition (\ref{h0}) 
explicitly, instead of hidden in the complementarity condition. This way of stating weak 
stationarity is also similar to that used in the MPCC setting, see 
\cite{AndreaniHaeserSecchinSilva,FlegelKanzow}.

In the next result we justify why Definition \ref{def:wstat_xy} is considered a weaker 
stationarity concept for the relaxed problem. 
\begin{theorem}
\label{th:kkt_wstat}
Suppose that $(\bar{x},\bar{y})$ is a KKT point for the relaxed problem (\ref{prob:relax1}). 
Then $(\bar{x},\bar{y})$ is $W_{I}$-stationary.
\end{theorem}
\beginproof
Denoting $\xi(x,y)=x*y$, we have 
$\nabla\xi_i(\bar{x},\bar{y})=\left(\begin{array}{c} \bar{y}_ie_i \\ 
\bar{x}_ie_i \end{array}\right)$. 
By the hypothesis, there exists a vector 
$$
(\lambda^g,\lambda^h,\lambda^{\theta},\mu,\lambda^{\tilde{H}},\lambda^\xi)\in
\R_+^m\times\R^p\times\R_+\times\R_+^n\times\R_+^n\times\R^n
$$ 
such that
\begin{align*}
\left(\begin{array}{c}\nabla f(\bar{x}) \\ 0 \end{array}\right)+
\sum_{i=1}^{m}\lambda_i^g\left(\begin{array}{c}\nabla g_i(\bar{x}) \\ 0
\end{array}\right)+\sum_{i=1}^{p}\lambda_i^h
\left(\begin{array}{c}\nabla h_i(\bar{x}) \\ 0 \end{array}\right)
+\lambda^{\theta}\nabla\theta(\bar{x},\bar{y})  \\
+\sum_{i=1}^{n}\mu_i \nabla H_i(\bar{x},\bar{y}) 
+\sum_{i=1}^{n}\lambda_i^{\tilde{H}}\nabla \tilde{H}_i(\bar{x},\bar{y})
+\sum_{i=1}^{n}\lambda_i^\xi\nabla\xi_i(\bar{x},\bar{y})
=\left(\begin{array}{c} 0 \\ 0 \end{array}\right),
\end{align*}
which means that 
\begin{align}
\nabla f(\bar{x})+\sum_{i=1}^{m}\lambda_i^g\nabla g_i(\bar{x})+\sum_{i=1}^{p}\lambda_i^h
\nabla h_i(\bar{x})+\sum_{i\in I_{01}}\lambda_i^\xi e_i
+\sum_{i\in I_{0+}}\lambda_i^\xi\bar{y}_ie_i=0, \label{nablaLrx} \\ 
-\lambda^{\theta}e-\mu+\lambda^{\tilde{H}}
+\sum_{i\in I_{\pm 0}}\lambda_i^\xi\bar{x}_ie_i=0, \label{nablaLry} 
\end{align}
where, for simplicity, we denoted $I_{0+}=I_{0+}(\bar{x},\bar{y})$ and 
$I_{\pm 0}=I_{\pm 0}(\bar{x},\bar{y})$. 

Moreover, we have 
\begin{align}
(\lambda^g)^Tg(\bar{x})=\lambda^{\theta}\theta(\bar{x},\bar{y})=
\mu^TH(\bar{x},\bar{y})=(\lambda^{\tilde{H}})^T\tilde{H}(\bar{x},\bar{y})=0. 
\label{compl_relax}
\end{align}
Defining 
$$
\begin{array}{rcc}
\lambda_i^G=\left\{\begin{array}{l} \lambda_i^\xi\bar{y}_i, \mbox{ for } 
i\in{I_{01}\cup I_{0+}}, \vspace{2pt}  \\ 
0, \mbox{ for } i\in {I_{\pm 0}\cup I_{00}}\end{array}\right.
&\mbox{and} &
\lambda_i^H=\left\{\begin{array}{l} \mu_i, \mbox{ for } i\in I_{0}, \vspace{2pt} \\ 
\mu_i-\lambda_i^\xi\bar{x}_i, \mbox{ for } i\in I_{\pm 0},\end{array}\right.
\end{array}
$$
we conclude immediately that $(\lambda^g,\lambda^h,\lambda^{\theta},
\lambda^H,\lambda^{\tilde{H}},\lambda^G)$ satisfies (\ref{nablaLx0}) and items 
(\ref{g0})--(\ref{h0}) of Definition \ref{def:wstat_xy}. To finish the proof, note that 
\begin{align*}
-\lambda^{\theta}e-\lambda^H+\lambda^{\tilde{H}} &  = 
-\lambda^{\theta}e+\lambda^{\tilde{H}}-\sum_{i\in I_{0}}\mu_i e_i-
\sum_{i\in I_{\pm 0}}(\mu_i-\lambda_i^\xi\bar{x}_i)e_i
\\ &  =
-\lambda^{\theta}e+\lambda^{\tilde{H}}-\mu+\sum_{i\in I_{\pm 0}}\lambda_i^\xi\bar{x}_ie_i
\end{align*}
which in view of (\ref{nablaLry}) gives (\ref{nablaLy0}). 
\endproof

As we have discussed before, a minimizer of the relaxed problem does not necessarily 
satisfy the KKT conditions mostly because of the complementarity constraint, which 
may prevent the fulfillment of constraint qualifications. This fact was 
illustrated in Example~\ref{ex:notGCQ}. Let us revisit this example in light of our 
$W_{I}$-stationarity concept. Now we can capture the minimizer by means of the KKT 
conditions for the tightened problem. 
\begin{example}
\label{ex:wstat}
Consider the MPCaC and the corresponding relaxed problem presented in Example \ref{ex:notGCQ}. 
$$
\begin{array}{lr}
\begin{array}{cl}
\displaystyle\mathop{\rm minimize }_{x\in\R^2} & x_1+x_2  \\
{\rm subject\ to } & -x_1+x_2^2\leq 0, \\
& \|x\|_0\leq 1, \\ & \\ &
\end{array}
&
\hspace{.5cm}
\begin{array}{cl}
\displaystyle\mathop{\rm minimize }_{x,y\in\R^2} & x_1+x_2 \\
{\rm subject\ to } & -x_1+x_2^2\leq 0, \\
& y_1+y_2\geq 1, \\
& x_iy_i=0,\; i=1,2, \\
& 0\leq y_i \leq 1, \; i=1,2.
\end{array}
\end{array}
$$
We saw that $x^*=(0,0)$ is the unique global solution of MPCaC and $(x^*,y^*)$, with  
$y^*=(1,0)$, is a global solution of the relaxed problem. Besides, this pair does not 
satisfy GCQ and it is not a KKT point. However, for $I=I_{0}(x^*)$, 
we conclude that $(x^*,y^*)$ is a KKT point for the tightened problem, that is, it is 
a $W_{I}$-stationary for the relaxed problem. Note also that ACQ holds at this 
point for TNLP$_{I}(x^*,y^*)$.
\end{example}

As a matter of fact, we can state $W_{I}$-stationarity using only the original 
variables $x$. This follows from the next result.
\begin{proposition}
\label{prop:wstat_x}
Let $(\bar{x},\bar{y})$ be feasible for the relaxed problem~(\ref{prob:relax1}). 
Then $(\bar{x},\bar{y})$ is $W_{I}$-stationary for this problem if and only if there exists a 
vector 
$$
(\lambda^g,\lambda^h,\gamma)\in\R_+^m\times\R^p\times\R^{|I|}
$$ 
such that 
\begin{align}
\nabla f(\bar{x})+\sum_{i=1}^{m}\lambda_i^g\nabla g_i(\bar{x})+
\sum_{i=1}^{p}\lambda_i^h\nabla h_i(\bar{x})+\sum_{i\in I}\gamma_ie_i=0, \label{gradx} \\ 
(\lambda^g)^Tg(\bar{x})=0. \label{compl_g}
\end{align}
\end{proposition}
\beginproof
Suppose first that $(\bar{x},\bar{y})$ is $W_{I}$-stationary, that is, it satisfies 
Definition \ref{def:wstat_xy}. In view of (\ref{nablaLx0}), if we define $\gamma=\lambda^G$, 
we obtain (\ref{gradx}). Moreover, (\ref{compl_g}) follows from 
Definition \ref{def:wstat_xy}(\ref{g0}).

Conversely, assume that $(\bar{x},\bar{y})$ satisfies conditions (\ref{gradx}) 
and (\ref{compl_g}). Then, defining $\lambda^G=\gamma$ and setting 
$\lambda^{\theta}=0$, $\lambda^H=\lambda^{\tilde{H}}=0$, we obtain (\ref{nablaLx0}) and 
(\ref{nablaLy0}). Therefore, $(\bar{x},\bar{y})$ is a $W_{I}$-stationary point.
\endproof

\begin{remark}
\label{rm:wstatx}
Note that the conditions (\ref{gradx}) and (\ref{compl_g}) generalize 
the concepts of $S$- and $M$-stationarity presented in \cite{BurdakovKanzowSchwartz16} 
if we consider $I=I_{0+}(\bar{x},\bar{y})\cup I_{01}(\bar{x},\bar{y})$ and 
$I=I_0(\bar{x})$, respectively. 
However, we stress that here we have different levels of weak stationarity, 
according to the set $I$ between $I_{0+}(\bar{x},\bar{y})\cup I_{01}(\bar{x},\bar{y})$ and 
$I_0(\bar{x})$. Moreover, an interesting and direct consequence of 
Proposition \ref{prop:wstat_x} is that stationarity gets stronger as the index set reduces 
(cf. proposition below).
\end{remark}

\begin{proposition}
\label{prop:wstat_order}
Let $(\bar{x},\bar{y})$ be feasible for the relaxed problem~(\ref{prob:relax1}). If 
$$
I_{0+}(\bar{x},\bar{y})\cup I_{01}(\bar{x},\bar{y})\subset I'\subset I\subset I_0(\bar{x}),
$$
then W$_{I'}$-stationarity implies $W_{I}$-stationarity. 
\end{proposition}

Proposition \ref{prop:wstat_x} also makes easier the task of verifying whether a 
point is or is not $W_{I}$-stationary, as we can see in the next example.
\begin{example}
\label{ex:wstatn}
Consider the following MPCaC and the associated relaxed problem. 
$$
\begin{array}{lr}
\begin{array}{cl}
\displaystyle\mathop{\rm minimize }_{x\in\R^n} & f(x)  \\
{\rm subject\ to } & g(x)\leq 0, \\
& \|x\|_0\leq n-1, \\ & \\ &
\end{array}
&
\hspace{.5cm}
\begin{array}{cl}
\displaystyle\mathop{\rm minimize }_{x,y\in\R^n} & f(x) \\
{\rm subject\ to } & g(x)\leq 0, \\
& e^Ty\geq 1, \\
& x*y=0, \\
& 0\leq y\leq e, 
\end{array}
\end{array}
$$
where the functions $f:\R^n\to\R$ and $g:\R^n\to\R^{n-1}$ are given by $f(x)=e^Tx$ 
and $g_i(x)=-x_i+x_n^2$, 
$i=1,\ldots,n-1$. Note that every feasible point $x$ of MPCaC satisfies $x_n=0$ because otherwise 
we would have $\|x\|_0=n$. Therefore, $x^*=0$ is the unique global solution of MPCaC and 
$(x^*,y^*)$, with $y^*=e_1$, is a global solution of the relaxed problem. Besides, we have 
$$
I_{0}=\{1,\ldots,n\},\;\; I_{01}=\{1\}, \;\;
I_{00}=\{2,\ldots,n\}\quad\mbox{and}\quad I_{0+}=\emptyset.
$$
Consider an arbitrary index set $I$ such that $I_{01}\subset I\subset I_{0}$. Then, the point 
$(x^*,y^*)$ is $W_{I}$-stationary if and only if $n\in I$.

Indeed, suppose first that $n\in I$.  Defining 
$\lambda_i^g=1$, $i=1,\ldots,n-1$, $\gamma_n=-1$ and $\gamma_i=0$, $i\in I\setminus\{n\}$, 
we have 
$$
\nabla f(x^*)+\sum_{i=1}^{n-1}\lambda_i^g\nabla g_i(x^*)+\sum_{i\in I}\gamma_ie_i=
\left(\begin{array}{c} \tilde e \\ 1 \end{array}\right)+
\left(\begin{array}{r} -\tilde e \\ 0 \end{array}\right)+
\left(\begin{array}{r} 0 \\ -1 \end{array}\right),
$$
where $\tilde e$ denotes the vector of all ones in $\R^{n-1}$. So, we obtain 
(\ref{gradx}) and (\ref{compl_g}). On the other hand, if $n\notin I$, then there is no 
vector $(\lambda^g,\gamma)\in\R_+^{n-1}\times\R^{|I|}$ such that the expression 
$$
\nabla f(x^*)+\sum_{i=1}^{n-1}\lambda_i^g\nabla g_i(x^*)+\sum_{i\in I}\gamma_ie_i=
\left(\begin{array}{c} \tilde e \\ 1 \end{array}\right)
-\left(\begin{array}{c} \lambda^g \\ 0 \end{array}\right)+
\left(\begin{array}{c} \gamma_{I} \\ 0 \end{array}\right)
$$
vanishes. Therefore, condition (\ref{gradx}) can never be satisfied. 
\end{example}

In view of Proposition \ref{prop:wstat_x}, we can relate our stationarity concept 
to other notions of stationarity, as $L$-, $N$- and $T$-stationarity proposed in 
\cite{BeckEldar,LiSong,PanXiuZhou}, for sparsity constrained optimization (when there 
is solely the cardinality constraint). Indeed, if $X=\mathbb{R}^n$, then a 
feasible point $(\bar{x},\bar{y})$ is $W_{I}$-stationary if and only if there exists 
a vector $\gamma\in\R^{|I|}$ such that 
$$
\nabla f(\bar{x})+\sum_{i\in I}\gamma_ie_i=0.
$$
Thus, $L$-stationarity \cite{BeckEldar} implies $W_{I_0}$-stationarity, which in turn is 
equivalent to $N^C$- and $T^C$-stationarity \cite{LiSong,PanXiuZhou}.

We conclude this section with an important consequence of the Remark \ref{rm:wstatx} 
and Proposition \ref{prop:wstat_order} that refers to the fulfillment of 
$W_{I}$-stationarity under some CQ for MPCaC. In \cite{CervinkaKanzowSchwartz} the 
authors introduced several MPCaC-tailored constraint qualifications (CC-CQ). 
In particular, they proved that $S$-stationarity (and 
consequently $M$-stationarity) \cite{BurdakovKanzowSchwartz16} holds at 
minimizers under each one of the proposed CC-CQ. Therefore, since $S$-stationarity 
at a point $(\bar{x},\bar{y})$ is equivalent to $W_{I_{S}}$-stationarity, with 
$I_S=I_{0+}(\bar{x},\bar{y})\cup I_{01}(\bar{x},\bar{y})$, and $W_{I_{S}}$ implies 
any other $W_{I}$ in the range $I_{S}\subset I\subset I_0(\bar{x})$, we have that 
a local minimizer of the relaxed problem (\ref{prob:relax}), for which a CC-CQ 
holds, satisfies $W_{I}$.

\subsubsection{Strong stationarity}
\label{sec:strong}
In this section we define strong stationarity as a special case of Definition \ref{def:wstat_xy} 
and prove, among other results, that it coincides with the classical KKT conditions for 
the relaxed problem. 
\begin{definition}
\label{def:sstat_xy}
Consider a feasible point $(\bar{x},\bar{y})$ for the relaxed problem (\ref{prob:relax1}). 
We say that $(\bar{x},\bar{y})$ is strongly stationary ($S$-stationary) for this problem if 
it is $W_{I}$-stationary, with $I=I_{0+}(\bar{x},\bar{y})\cup I_{01}(\bar{x},\bar{y})$,  
according to Definition \ref{def:wstat_xy}.
\end{definition}

Let us start by showing the equivalence between $S$-stationarity and KKT. 
\begin{theorem}
\label{th:sstat_kkt}
Let $(\bar{x},\bar{y})$ be a feasible point for the relaxed problem (\ref{prob:relax1}). Then, 
$(\bar{x},\bar{y})$ is $S$-stationary if and only if it satisfies the usual KKT conditions for 
this problem.
\end{theorem}
\beginproof
Note first that sufficiency follows directly from Theorem \ref{th:kkt_wstat}. 
To prove the necessity, note that by Proposition \ref{prop:wstat_x} there exists a vector 
$$
(\lambda^g,\lambda^h,\gamma)\in\R_+^m\times\R^p\times\R^{|I_{0+}\cup I_{01}|}
$$ 
such that 
\begin{align*}
\nabla f(\bar{x})+\sum_{i=1}^{m}\lambda_i^g\nabla g_i(\bar{x})+
\sum_{i=1}^{p}\lambda_i^h\nabla h_i(\bar{x})+\sum_{i\in I_{0+}\cup I_{01}}\gamma_ie_i=0, \\ 
(\lambda^g)^Tg(\bar{x})=0.
\end{align*}
Defining $\lambda^{\theta}=0$, $\lambda^{H}=\lambda^{\tilde{H}}=0$ and 
$$
\lambda_i^\xi=\left\{\begin{array}{l} \gamma_i/\bar{y}_i, \mbox{ for } 
i\in I_{0+}\cup I_{01}, \vspace{2pt} \\ 
0, \mbox{ for } i\in I_{\pm 0}\cup I_{00}, \end{array}\right.
$$
we conclude that the vector 
$$
(\lambda^g,\lambda^h,\lambda^{\theta},\lambda^H,\lambda^{\tilde{H}},\lambda^\xi)\in
\R_+^m\times\R^p\times\R_+\times\R_+^n\times\R_+^n\times\R^n
$$ 
fulfills the conditions 
\begin{align}
\nabla f(\bar{x})+\sum_{i=1}^{m}\lambda_i^g\nabla g_i(\bar{x})+\sum_{i=1}^{p}\lambda_i^h
\nabla h_i(\bar{x})+\sum_{i\in I_{0+}\cup I_{01}}\lambda_i^\xi\bar{y}_ie_i=0, 
\label{gradx_relax1} \\ 
-\lambda^{\theta}e-\lambda^H+\lambda^{\tilde{H}}
+\sum_{i\in I_{\pm 0}}\lambda_i^\xi\bar{x}_ie_i=0, \label{grady_relax1} \\
(\lambda^g)^Tg(\bar{x})=\lambda^{\theta}\theta(\bar{x},\bar{y})=
(\lambda^H)^TH(\bar{x},\bar{y})=
(\lambda^{\tilde{H}})^T\tilde{H}(\bar{x},\bar{y})=0. \label{compl_relax1}
\end{align}
Since we have defined $\xi(x,y)=x*y$, we see that (\ref{gradx_relax1}) and (\ref{grady_relax1}) 
are equivalent to 
\begin{align*}
\left(\begin{array}{c}\nabla f(\bar{x}) \\ 0 \end{array}\right)+
\sum_{i=1}^{m}\lambda_i^g\left(\begin{array}{c}\nabla g_i(\bar{x}) \\ 0
\end{array}\right)+\sum_{i=1}^{p}\lambda_i^h
\left(\begin{array}{c}\nabla h_i(\bar{x}) \\ 0 \end{array}\right)
+\lambda^{\theta}\nabla\theta(\bar{x},\bar{y})  \\
+\sum_{i=1}^{n}\lambda_i^H \nabla H_i(\bar{x},\bar{y}) 
+\sum_{i=1}^{n}\lambda_i^{\tilde{H}}\nabla \tilde{H}_i(\bar{x},\bar{y})
+\sum_{i=1}^{n}\lambda_i^\xi\nabla\xi_i(\bar{x},\bar{y})
=\left(\begin{array}{c} 0 \\ 0 \end{array}\right)
\end{align*}
and, consequently, $(\bar{x},\bar{y})$ is a KKT point for the relaxed problem (\ref{prob:relax1}). 
\endproof

As an immediate consequence of Theorem \ref{th:sstat_kkt} and Proposition \ref{prop:wstat_kkttnlp} 
we have the equivalence between the KKT conditions for the relaxed and tightened problems when 
$I=I_{0+}(\bar{x},\bar{y})\cup I_{01}(\bar{x},\bar{y})$.
\begin{corollary}
\label{c:kkt_kkt}
Consider a feasible point $(\bar{x},\bar{y})$ for the relaxed problem (\ref{prob:relax1}) and 
$I=I_{0+}(\bar{x},\bar{y})\cup I_{01}(\bar{x},\bar{y})$. Then, $(\bar{x},\bar{y})$ is a 
KKT point for (\ref{prob:relax1}) if and only if it is a KKT point for the tightened 
problem TNLP$_{I}(\bar{x},\bar{y})$, given by (\ref{prob:tight}). 
\end{corollary}

It should be noted that, if $I_{00}(\bar{x},\bar{y})=\emptyset$ then all stationarity conditions 
presented here are the same and correspond to KKT for the relaxed problem. This follows directly 
from (\ref{indexI}). In this case, we say that $(\bar{x},\bar{y})$ satisfies the 
{\em strict complementarity}.

On the other hand, given an arbitrary index set $I$ satisfying (\ref{indexI}), we have the 
following result.
\begin{proposition}
\label{prop:lbdGI00_kkt}
Suppose that $(\bar{x},\bar{y})$ is $W_{I}$-stationary for the relaxed 
problem (\ref{prob:relax1}) with the associated vector of multipliers 
$\lambda=(\lambda^g,\lambda^h,\lambda^{\theta},\lambda^H,\lambda^{\tilde{H}},\lambda^G)$. 
If $\lambda_{I_{00}(\bar{x},\bar{y})}^G=0$, then $(\bar{x},\bar{y})$ is $S$-stationary. 
\end{proposition}
\beginproof
Note that Definition \ref{def:wstat_xy}, relation (\ref{nablaLx0}) and the condition 
$\lambda_{I_{00}(\bar{x},\bar{y})}^G=0$ allows us to write 
\begin{align*}
\nabla f(\bar{x})+\sum_{i=1}^{m}\lambda_i^g\nabla g_i(\bar{x})+
\sum_{i=1}^{p}\lambda_i^h\nabla h_i(\bar{x})+\sum_{i\in I_{0+}\cup I_{01}}\lambda_i^Ge_i=0, \\ 
(\lambda^g)^Tg(\bar{x})=0.
\end{align*}
Using Proposition \ref{prop:wstat_x} with $I=I_{0+}\cup I_{01}$, we see that 
$(\bar{x},\bar{y})$ is $W_{I}$-stationary, which means 
$S$-stationary. 
\endproof

Clearly (see Proposition \ref{prop:wstat_order}) $S$-stationarity implies W-stationarity. 
Moreover, this implication is strict in view of Example \ref{ex:wstat}.

\section{Conclusion and future works}
\label{sec:concl} 
In this paper we have presented a new and weaker stationarity condition, called 
$W_I$-stationarity, for Mathematical Programs with Cardinality Constraints (MPCaC). 
We proposed a unified approach that goes from the weakest to the strongest stationarity. 
Several theoretical results were presented, such as properties and relations 
concerning the reformulated problems and the original one. 
In addition, we discussed the relaxed problem by analyzing the classical constraints 
in two cases, linear and nonlinear, with results, examples and counterexamples. 

As we have seen in this paper, KKT implies $W_{I}$-stationarity. However, even being 
weaker, it is not a necessary optimality condition. Therefore, a subject of ongoing 
research is the proposal of sequential optimality conditions for MPCaC problems in 
the hope of obtaining a weaker condition than $W_{I}$-stationarity, which 
would be satisfied at every minimizer, independently of any constraint 
qualification. This in turn would allow us to discuss algorithmic consequences, once 
we believe that MPCaC is an important problem, both theoretically and numerically.

\vspace{.75cm}

\noindent{\bf Acknowledgements}
We thank the anonymous referees who helped us to improve the presentation of this work.



\end{document}